\documentclass[twocolumn,twoside]{IEEEtran}
\usepackage{amsmath,amssymb,amsthm,epsfig,dsfont,color,empheq,graphicx}
\usepackage{enumerate,url,algpseudocode,algorithm,wasysym,epstopdf}
\usepackage[table]{xcolor}
\usepackage{accents}
\usepackage[utf8]{inputenc}
\newcommand{\ubar}[1]{\underaccent{\bar}{#1}}
\usepackage{subfig}

\DeclareMathOperator{\diag}{diag}

\DeclareMathOperator{\real}{Re}

\newtheorem{proposition}{Proposition}

\theoremstyle{remark}

\begin{document}
\title{Scalable Electric Vehicle Charging Protocols}

\author{
	Liang Zhang,~\IEEEmembership{Student Member,~IEEE,}
	Vassilis Kekatos,~\IEEEmembership{Member,~IEEE,}
	and Georgios B. Giannakis,~\IEEEmembership{Fellow,~IEEE}

\vspace*{-0.7em}
	
\thanks{Work in this paper was supported by NSF grants 1423316, 1442686, 1508993, and 1509040. L. Zhang and G. B. Giannakis are with the Digital Technology Center and the ECE Dept., University of Minnesota, Minneapolis, MN 55455, USA.  V. Kekatos is with the ECE Dept., Virginia Tech, Blacksburg, VA 24061, USA. Emails: zhan3523@umn.edu, kekatos@vt.edu, georgios@umn.edu.}}

\maketitle

\begin{abstract}
Although electric vehicles are considered a viable solution to reduce greenhouse gas emissions, their uncoordinated charging could have adverse effects on power system operation. Nevertheless, the task of optimal electric vehicle charging scales unfavorably with the fleet size and the number of control periods, especially when distribution grid limitations are enforced. To this end, vehicle charging is first tackled using the recently revived Frank-Wolfe method. The novel decentralized charging protocol has minimal computational requirements from vehicle controllers, enjoys provable acceleration over existing alternatives, enhances the security of the pricing mechanism against data attacks, and protects user privacy. To comply with voltage limits, a network-constrained EV charging problem is subsequently formulated. Leveraging a linearized model for unbalanced distribution grids, the goal is to minimize the power supply cost while respecting critical voltage regulation and substation capacity limitations. Optimizing variables across grid nodes is accomplished by exchanging information only between neighboring buses via the alternating direction method of multipliers. Numerical tests corroborate the optimality and efficiency of the novel schemes.
\end{abstract}


\section{Introduction}\label{sec:intro}
Electric vehicles (EVs) have received significant attention from the automotive industry and the government due to their capacity to reduce greenhouse gas emissions and mitigate oil dependency. Nevertheless, the overall load profile will be greatly affected with increasing numbers of EVs. Uncoordinated charging of even a 10\% penetration of EV loads will notably affect power system operation, giving rise to voltage magnitude fluctuations and unacceptable load peaks~\cite{RFK12}. On the other hand, with proper coordination scheme, EV loads can be controlled to minimize charging costs or perform valley-filling tasks relying on advanced power electronics.

Different charging control schemes have been proposed. A centralized scheduling scheme to minimize total charging costs based on the time-of-use price has been devised in \cite{YJCao12}. However, new load peaks may arise during low-price (also termed valley) periods. In \cite{SPR09}, vehicle plug-in times are decided using random numbers, hence neglecting the specific charging requests of individual EV users. Charging rates have been also optimized in a centralized manner to facilitate voltage regulation~\cite{RFK12}, yet the number of control variables scales unfavorably with the number of vehicles.

Decentralized control strategies not only offer computational and communication savings, but they oftentimes enhance the privacy of vehicle users since they do not require the charging requests of EVs to be uploaded to the control center. Decentralized charging protocols are available based on congestion pricing schemes similar to those used in Internet Protocol (IP) networks; nevertheless, their optimality is not guaranteed~\cite{IP13}. Presuming identical plug-in/-out times and energy requests for all vehicles, a game-theoretic charging scheme attaining a Nash equilibrium has been developed in~\cite{mah}. Iterative schemes based on Lagrangian relaxation are suggested in~\cite{NGG}, while \cite{Rivera13} builds on the alternating direction method of multipliers (ADMM). Distribution locational marginal prices are leveraged to coordinate vehicle charging in~\cite{RQO14}. Reference~\cite{GTL11} proves a feasible valley-filling charging profile to be optimal for any convex cost, and it develops a decentralized protocol. A multi-agent based EV charging scheme is proposed in~\cite{karfopoulos2013multi}. 

Vehicle charging under distribution grid limitations has been studied too. Centralized EV scheduling is studied under different linear models for multiphase networks in~\cite{NCEV} and~\cite{MIPEV}. The objective function is confined to be linear and the optimal solution is found using generic commercial solvers without exploiting the problem structure. A method for heuristically checking network constraint violations after vehicles have been scheduled is reported in~\cite{rev}. Presuming at most one EV per bus, management under balanced network constraints has been tackled using a water-filling algorithm~\cite{NCQ14}.

The optimal vehicle charging problem considered here can be rigorously stated as follows. Given charging requests from EVs across time, a utility company schedules their charging to minimize certain cost function, e.g., the power supply cost or the load variance. The latter is equivalent to the so termed the \emph{valley-filling} task. Depending on whether grid specifications are taken into account, two charging scenarios can be recognized. The first scenario ignores any grid-related constraints. Such a scenario arises for example when the EV load is relatively low and is not expected to incur voltage or feeder violations; see the valley-filling task in~\cite{GTL11}. In this first scenario, vehicle charging may be alternatively performed by a charging station or a load aggregator to minimize its power supply cost. Under the second scenario, EV penetration is high, and thus, voltage regulation and feeder limitations must be enforced by the utility. Apparently, the first scenario constitutes a relaxation of the second scenario of \emph{network-constrained} vehicle charging. Thus, protocols for the first scenario will be used as building modules for the second one. 

Our contribution is two-fold. First, a decentralized charging method based on the Frank-Wolf algorithm is developed (Section~\ref{sec:singlebus}). Different from existing schemes, the novel protocol requires minimal requirements from the vehicle charging controllers and involves privacy-preserving updates. Numerical tests demonstrate that the closed-form low-complexity updates yield significant convergence improvement over existing alternatives (Section~\ref{subsec:simua1}). Secondly, building on an approximate distribution grid model, network-constrained EV charging is formulated as a convex quadratic program (Section~\ref{sec:network}), and tackled using a decentralized scheme based on ADMM (Section~\ref{sec:jo-admm}). Compared to existing centralized schemes, the decentralized protocol requires communication only between neighbors and preserves the privacy of EV owners. Numerical tests on the unbalanced IEEE 123-bus feeder corroborate the optimality of the proposed charging protocol (Section~\ref{subsec:simua1ncev}). 

Regarding \emph{notation}, column vectors (matrices) are denoted by lower- (upper-) case boldface letters, except for power flow vectors $(\mathbf{P},\mathbf{Q},\mathbf{S})$. Sets are represented using calligraphic symbols, and $|\mathcal{S}|$ is the cardinality of $\mathcal{S}$. Symbol $^{\top}$ stands for transposition; while $\mathbf{0}$, $\mathbf{1}$, and $\mathbf{e}_n$, denote respectively the all-zeros, all-ones, and the $n$-th canonical vectors. Operator $\diag(\mathbf{x})$ defines a diagonal matrix having $\mathbf{x}$ on its diagonal, and $\real(z)$ returns the real part of complex number $z$.

\section{Optimal Vehicle Charging}\label{sec:singlebus}
This section studies EV charging without network constraints. Under this scenario, the utility company, a load aggregator, or a charging station would like to coordinate EVs to minimize the power supply cost or for valley-filling purposes. Upon formulating the problem, an optimal charging scheme is developed and contrasted to state-of-the-art alternatives.
  
\subsection{Electric Vehicle Charging Model}\label{subsec:ev-model}
An EV scheduler
coordinates the charging of $M$ EVs over a period of $T$ consecutive time slots comprising the set $\mathcal{T}:=\{1,\ldots,T\}$. The time slot duration $\Delta T$ can range from minutes to an hour, depending on charging specifications, the granularity of load forecasts, as well as communication and computation capabilities. Let $e_m(t)$ denote the energy charge for vehicle $m$ at time $t$ with $m=1,\ldots,M$, and $t\in\mathcal{T}$. Given that operational slots have equal duration, the terms power and energy will be used interchangeably. The charge $e_m(t)$ can range from zero to its maximum value $\bar{e}_m(t)$. Apparently, a vehicle can be charged only when it is connected to the grid. If $\mathcal{T}_m\subseteq \mathcal{T}$ is the set of time slots that vehicle $m$ is connected to the grid (not necessarily consecutive), then for all $t\in\mathcal{T}$
 \begin{equation*}
 \bar{e}_m(t)=\left\{
 \begin{array}{ll}
 \bar{e}_m&,~t\in \mathcal{T}_m\\
 0 &,~ \text{otherwise}
 \end{array} \right.
 \end{equation*}
 where $\bar{e}_m$ is the maximum charging rate determined by the battery of vehicle $m$. Let $\mathbf{e}_m:=[e_m(1) ~ \cdots ~ e_m(T)]^{\top}$ be the charging profile for EV $m$. Profile $\mathbf{e}_m$ should belong to the compact and convex set
 \begin{equation}\label{eq:evc}
 \mathcal{E}_m:=\{\mathbf{e}_m:\mathbf{e}_m^{\top}\mathbf{1}=R_m,~0 \leq e_m(t) \leq \bar{e}_m(t)~\forall~ t\in\mathcal{T}\}
 \end{equation} 
where $R_m$ is the total energy needed by EV $m$. The latter depends on the initial state of charge, the desired state of charge, and the efficiency of the battery.

Through coordinated charging of electric vehicles, various objectives can be achieved, such as minimizing charging costs or valley-filling. Optimal EV charging can be posed as the optimization problem~\cite{GTL11}
\begin{align}\label{eq:EVS2}
\min_{\{\mathbf{e}_m\}_{m=1}^M} ~&~C(\{\mathbf{e}_m\}):=\sum_{t=1}^T C_t\left(d(t)+\sum_{m=1}^M e_m(t)\right)\\
\textrm{s.to}~&~\mathbf{e}_m \in\mathcal{E}_m,~\forall~m=1,\ldots,M\nonumber
\end{align}
where the energy costs $C_t(\cdot): \mathbb{R}\rightarrow\mathbb{R}$ are assumed convex and differentiable. For charging cost minimization, $\{C_t\}_{t=1}^T$ can be linear or quadratic~\cite{YJCao12}; while $C_t(x)=x^2/2$ for all $t$ when it comes to the valley-filling task. 
Parameters $\{d(t)\}^T_{t=1}$ capture the based load for the EV scheduler, which is assumed inelastic and known in advance. The network constrained EV charging is postponed for Section~\ref{sec:network}, wherein problem \eqref{eq:EVS2} turns out to be a building module. To facilitate scheduling, each electric vehicle controller is capable of two-way communication and of executing simple computation tasks. Before the beginning of the charging horizon $\mathcal{T}$, vehicle controller submit their charging requests $\{(\mathcal{T}_m,R_m)\}$ to the charging station controller. Protocols for efficiently solving \eqref{eq:EVS2} are presented next.
 
\subsection{Scalable Charging Protocol}\label{subsec:ev-fw}
Observe that the total number of variables involved in \eqref{eq:EVS2} is $MT$. Therefore, although \eqref{eq:EVS2} is a convex problem, solving it is a non-trivial task, particularly for large EV fleets and/or decreasing control intervals $\Delta T$. To derive a scalable solver, the Frank-Wolfe method is deployed next~\cite{jaggi2013revisiting}. Also known as conditional gradient algorithm, the Frank-Wolfe method aims at solving the generic problem
 \begin{equation}\label{eq:problem}
 \mathbf{y}^*\in\arg\min_{\mathbf{y}\in \mathcal{Y}}~f(\mathbf{y})
 \end{equation}
 where $f$ is a differentiable convex function, and $\mathcal{Y}$ is a compact convex set. The method selects an initial $\mathbf{y}^0\in\mathcal{Y}$, and iterates between the updates for $k=1,2,\ldots,$ as
 \begin{subequations}\label{eq:FW}
 \begin{align}
 \mathbf{r}^k&\in\arg\min_{\mathbf{r} \in \mathcal{Y}} ~ \mathbf{r}^{\top} \nabla f(\mathbf{y}^k)\label{eq:FW:g}\\
 \mathbf{y}^{k+1}&:=\mathbf{y}^k+\eta_k(\mathbf{r}^k-\mathbf{y}^k)\label{eq:FW:x}
 \end{align}
 \end{subequations}
 with $\eta_k:=2/(k+2)$. Step \eqref{eq:FW:g} finds $\mathbf{r}^k$ such that $(\mathbf{r}^k-\mathbf{y}^k)$ is a feasible descent direction for the first-order approximation of the cost in \eqref{eq:problem}. Step \eqref{eq:FW:x} updates $\mathbf{y}^k$ towards that direction after scaling it with the diminishing step size $\eta_k$. The updated $\mathbf{y}^{k+1}$ is always feasible, since it is computed as the convex combination of $\mathbf{y}^{k}\in\mathcal{Y}$ and $\mathbf{r}^k\in\mathcal{Y}$.
 
Granted that \eqref{eq:EVS2} entails a differentiable cost and a compact feasible set; it is amenable to Frank-Wolfe iterations. In the first Frank-Wolfe step, the gradient of the cost in \eqref{eq:EVS2} with respect to $\{\mathbf{e}_m\}_{m=1}^M$ must be obtained. Critically, due to the problem structure, the per-vehicle partial gradients of the cost are all \emph{identical} to
 \begin{equation*}
 \nabla_{\mathbf{e}_m}C(\{\mathbf{e}_m\})=\mathbf{g},\quad m=1,\ldots,M.
 \end{equation*}
 It can be readily checked that the $t$-th entry of the common partial gradient $\mathbf{g}\in\mathbb{R}^T$ evaluated at $\{\mathbf{e}_m^k\}$ is
 \begin{equation}\label{eq:grad2}
 g^k(t)=\nabla_{e_m^k(t)} C_t \left(d(t)+\sum_{m=1}^M e_m^k(t)\right),~t=1,\ldots,T.
 \end{equation}
 Applying \eqref{eq:FW:g} to the problem at hand requires solving
 \begin{equation}\label{eq:lp}
 \{\mathbf{r}_m^k\}_{m=1}^M\in\arg\min_{\{\mathbf{r}_m\in \mathcal{E}_m\}_{m=1}^M}~\sum_{m=1}^M \mathbf{r}_m^{\top} \mathbf{g}^k
 \end{equation}
 which is separable across vehicles. Thus, given $\mathbf{g}^k$, vehicle $m$ needs to solve the linear program
 \begin{equation}\label{eq:lp2}
 \mathbf{r}_m^k\in\arg\min_{\mathbf{r}_m\in \mathcal{E}_m}~\mathbf{r}_m^{\top} \mathbf{g}^k.
 \end{equation}
Problem \eqref{eq:lp2} involves a linear cost minimized over a weighted budget and box constraints. The key observation here is that due to the aforementioned structure, problem \eqref{eq:lp2} can be solved by a simple sorting algorithm~\cite[Chap.~4]{BoVa04}: The entries of $\mathbf{g}^k$ are first sorted in increasing order as 
 \begin{equation} \label{eq:torder}
 g^k(t_1^k)\leq g^k(t_2^k)\leq \ldots \leq g^k(t_T^k).
 \end{equation}
Since the problems in \eqref{eq:lp2} share vector $\mathbf{g}^k$ for all $m$, the sorting operation is performed only once by the charging station. Then, for vehicle $m$, we need to find the index $J_m^k$ for which
 \begin{equation}\label{eq:opts0}
 \sum_{j=1}^{J_m^k}\bar{e}_m(t_j^k)\leq R_m~\text{and}~\sum_{j=1}^{J_m^k+1}\bar{e}_m(t_j^k)> R_m.
 \end{equation}
 Subsequently, the entries of the minimizer $\mathbf{r}_m^k$ of \eqref{eq:lp2} can be computed per vehicle $m$ as
 \begin{equation}\label{eq:opts}
 r_m^{k}(t_j^k)=\left\{
 \begin{array}{ll}
 \bar{e}_m(t_j^k)&,~j=1,\ldots,J_m^k-1\\
 R_m-\sum_{j=1}^{J_m^k-1}\bar{e}_m(t_j^k)&,~j=J_m^k\\
 0&,~j=J_m^k+1,\ldots,T
 \end{array}\right..
 \end{equation}
 The solution in \eqref{eq:opts} simply selects the maximum possible charge during the cheapest time slots in a greedy fashion. Interestingly, finding $\mathbf{r}_m^k$ from \eqref{eq:opts} requires knowing solely the rank order (smallest to largest) rather than the actual entries of the gradient vector $\mathbf{g}^k$.

\begin{figure}[t]
\centering
\subfloat[Control center broadcasts time slot pricing ordering (from cheapest to most expensive) to EV controllers.]{\includegraphics[width=0.240\textwidth]{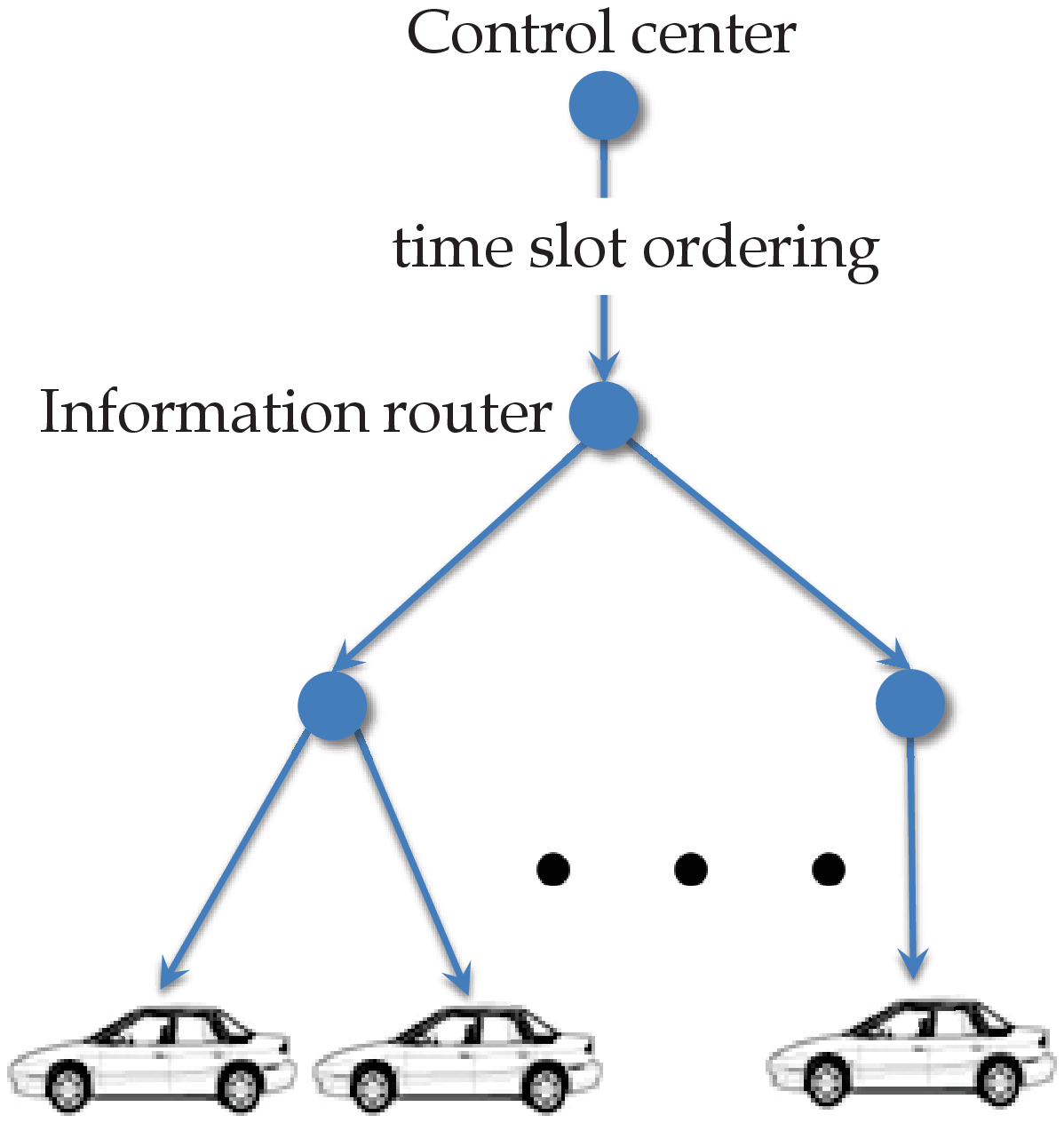}\label{fig:subf1}} ~ 
\subfloat[Summations of charging profiles are transmitted to charging control center.]{\includegraphics[width=0.258\textwidth]{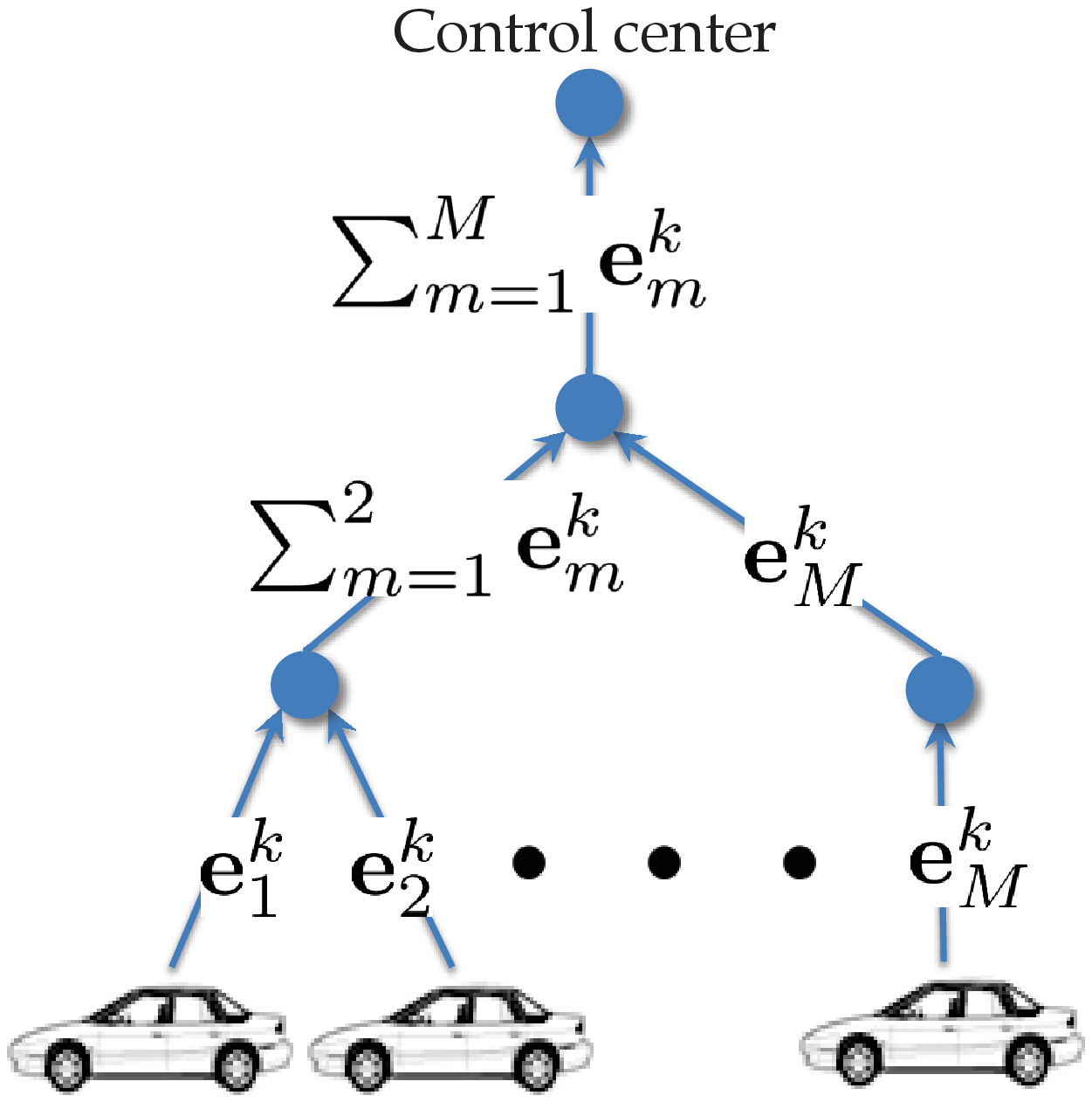}\label{fig:subf2}}
\caption{Information exchange for Algorithm~\ref{alg:docsev} at iteration $k$.}\label{fig:eva}
\end{figure}

The second Frank-Wolfe step updates the charging profiles via the convex combinations
 \begin{equation}\label{eq:cc}
 \mathbf{e}_m^{k+1}=(1-\eta_k)\mathbf{e}_m^k+\eta_k\mathbf{r}_m^k
 \end{equation}
for all vehicles $m=1,\ldots,M$.

 \begin{algorithm}[t]
 \caption{Decentralized EV scheduling}\label{alg:docsev}
 \begin{algorithmic}[1]
 \renewcommand{\algorithmicrequire}{\textbf{Input:}}
 \renewcommand{\algorithmicensure}{\textbf{Output:}} 
 \State Initialize $\mathbf{e}^0=\mathbf{0}$ and $\mathbf{g}^0(t)$ from~\eqref{eq:grad2}.
 \For{$k = 1,2,\ldots$}
 \State EV scheduler calculates $\mathbf{g}^k$ from \eqref{eq:grad2}.
 \State EV scheduler broadcasts $\mathbf{g}^k$ entry ranking to EVs.
 \State Vehicles update $\{\mathbf{e}_m^{k}\}_{m=1}^M$ via \eqref{eq:opts0}--\eqref{eq:cc}.
 \State Profile sums $\sum_{m=1}^{M}\mathbf{e}_m^{k}$ sent to control center. 
 \EndFor
 \end{algorithmic}
 \end{algorithm}

To practically implement \eqref{eq:grad2}--\eqref{eq:cc} during iteration $k$, the charging control center evaluates the cost gradients $\{g^k(t)\}_{t\in \mathcal{T}}$ defined in \eqref{eq:grad2}, and sorts them to determine the time slot ordering $\{t_1^k,t_2^k,\ldots,t_T^k\}$. This sorting operation can be performed using for example the Merge-Sort algorithm in $\mathcal{O}(T\log T)$ operations~\cite{donald1999art}. The price ordering of time slots is subsequently broadcast to all EV controllers as shown in Fig.~\ref{fig:subf1}. Based on its charging needs $\mathcal{E}_m$, the $m$-th EV controller first finds $\mathbf{r}_m^k$ from \eqref{eq:opts0}--\eqref{eq:opts} in $\mathcal{O}(T)$. It then updates its charging profile $\mathbf{e}_m^{k+1}$ using \eqref{eq:cc} in $\mathcal{O}(T)$. Note that operations \eqref{eq:opts0}--\eqref{eq:cc} can be performed in parallel over the $M$ EV controllers. The updated charging profiles $\{\mathbf{e}_m^{k+1}\}_{m=1}^M$ are communicated back to the charging center, where upon adding the base load $\{d(t)\}$, the center computes the updated cost gradient $\mathbf{g}^{t+1}$, and iterations proceed as tabulated in Algorithm~\ref{alg:docsev}. The developed solver converges to optimal charging profiles $\{\mathbf{e}_m^*\}$ at the rate~\cite{jaggi2013revisiting}
\begin{equation}
C(\{\mathbf{e}_m^k\})-C(\{\mathbf{e}_m^*\})\leq\mathcal{O}\left(\frac{1}{k}\right).
\end{equation}

Algorithm~\ref{alg:docsev} not only exhibits provable convergence and low computational cost (namely $\mathcal{O}(T\log T)$ operations) per iteration. It further enjoys two additional advantages. First, the charging center does not require knowing the individual charging profiles $\{\mathbf{e}_m^k\}$, since their summation $\sum_{m=1}^M \mathbf{e}_m^k$ suffices for finding the gradient vector $\mathbf{g}^k$. In an effort to preserve the privacy of EV users, a simple communication protocol can be designed. Information flow can be arranged over a tree graph rooted at the charging center, and vehicle controllers constitute the remaining tree nodes. Each node receives aggregate charging profiles from its downstream nodes, adds them up to its own profile, and forwards the updated aggregate charging profile to its parent node. As a second feature, vehicle controllers do not need to know the precise value of the cost gradient vector $\mathbf{g}^k$, but only the ordering of its entries (current price ordering of time slots). This algorithmic feature lightens the communication load from the charging center to the vehicles, and enhances resiliency to price manipulations and data attacks to the solving scheme.

 
\subsection{Comparison with Previous Work}\label{subsec:FWvsPGD}
The optimal EV charging of \eqref{eq:EVS2} has been previously studied in~\cite{GTL11}, where a projected gradient descent (PGD) solver was developed. Interpreted here as a projected gradient algorithm applied to minimize the non-strongly convex cost in \eqref{eq:EVS2}, the PGD method exhibits a convergence rate of $\mathcal{O}(\frac{1}{k})$~\cite{Be15}. At iteration $k$ of the PGD method, controller $m$ solves in parallel
  \begin{equation}\label{eq:PGD}
  \mathbf{e}_m^{k+1}:=\arg\min_{\mathbf{e}_m\in\mathcal{E}_m}~\|\mathbf{e}_m-(\mathbf{e}_m^k-\eta_k'\mathbf{g}^k) \|_2^2
  \end{equation}
for a step size $\eta_k'>0$. In other words, every EV controller projects vector $(\mathbf{e}_m^k-\eta_k'\mathbf{g}^k)$ onto the simplex $\mathcal{E}_m$, which is a non-trivial task.
  
On the other hand, each iteration of Algorithm~\ref{alg:docsev} involves closed-form updates, offering high computational efficiency and posing affordable hardware requirements on EV controllers. Although both Algorithm~\ref{alg:docsev} and the PGD solver are decentralized schemes with convergence rate $\mathcal{O}(\tfrac{1}{k})$, the overall computation time for the former is significantly lower due to its simpler per-iteration updates: The numerical tests in Section~\ref{subsec:simua1} demonstrate that Alg.~\ref{alg:docsev} provides a 100 times speed-up advantage over the PGD solver and the centralized solver SeDuMi. The SeDuMi solver would be a viable option for tackling \eqref{eq:EVS2} in a centralized manner after collecting all charging needs $\{(R_m,\mathcal{E}_m)\}$ at the charging center.

\section{Network-Constrained EV Scheduling}\label{sec:network}
The charging scheme of Section~\ref{sec:singlebus} applies to scenarios where EV charging can be transparently supported by the underlying grid. If higher levels of EV load incur voltage magnitude or feeder capacity violations, the underlying power distribution grid needs to be taken into account. In this context, upon reviewing an approximate model for unbalanced distribution grids, this section formulates a network-constrained vehicle charging task, while a decentralized solver scalable to the number of buses and EVs is developed in Section~\ref{sec:jo-admm}.

\subsection{Modeling Unbalanced Distribution Grids}\label{subsec:mp-model}
Electric vehicles are connected to a distribution feeder comprising $N+1$ buses indicated by $n\in\mathcal{N}:=\{0,1,\ldots,N\}$, and phases indexed by $\phi\in\{a,b,c,\}$. Let $\mathcal{M}_{n,\phi}$ represent the set of vehicles located on phase $\phi$ of bus $n$, and $M_{n,\phi}:=|\mathcal{M}_{n,\phi}|$. The distribution grid is assumed to be functionally radial with the substation bus numbered by $n=0$. Every non-feeder bus $n\in\mathcal{N}^+$ with $\mathcal{N}^+:=\mathcal{N}\setminus\{0\}$ has a unique parent bus indexed by $\pi_n$. The distribution line connecting bus $\pi_n$ with bus $n$ is denoted by $n$. For bus $n$, let also $\mathcal{C}_n$ denote the set of its children buses, and $\mathcal{P}_n \subseteq \{a,b,c\}$ the set of its phases.

To enforce distribution network and voltage regulation limitations, the underlying physical system is taken into account. For that purpose, the distribution grid can be captured either by the full AC power flow model or the linearized power flow model proposed in~\cite{LAGan}. The former becomes tractable under appropriate conditions using convex relaxations~\cite{HaozhuD}, \cite{LAGan}. However, counterexamples indicate that convex relaxations are not always successful and they can increase computational requirements. On the other hand, several numerical tests indicate that the approximation error of the linearized model is within the order of $10^{-3}$ in terms of calculating voltages~\cite{LAGan}, \cite{VZ}. Although the linearized grid model is adopted here to simplify calculations, extending our charging protocol to the full AC model is straightforward.

Let $\mathbf{v}_n$,  $\mathbf{p}_n$, and $\mathbf{q}_n$ be respectively the 3-dimensional vectors of squared voltage magnitudes and (re)active power injections for all phases of bus $n$. For line $n\in\mathcal{N}$, let $\mathbf{Z}_n=\mathbf{Z}_n^{\top}\in\mathbb{R}^{3\times 3}$ be the related phase impedance matrix, and $\mathbf{P}_n$ and $\mathbf{Q}_n$ be the vectors of (re)active power flows on all phases of line $n$. 
If line losses are relatively small and voltages are roughly balanced, the linearized multi-phase power flow model reads~\cite{LAGan}, \cite{VZ}
\begin{subequations}\label{eq:m-model}
\begin{align}
\mathbf{p}_n&=\sum_{k\in\mathcal{C}_n} \mathbf{P}_k  - \mathbf{P}_n\label{eq:m-amodelp}\\
\mathbf{q}_n&=\sum_{k\in\mathcal{C}_n} \mathbf{Q}_k  - \mathbf{Q}_n\label{eq:m-amodelq}\\
\mathbf{v}_{\pi_n} - \mathbf{v}_n &=  \real\left\{\bar{\mathbf{Z}}_n(\mathbf{P}_n+j\mathbf{Q}_n)\right\} \label{eq:m-amodelv}
\end{align}
\end{subequations}
where $\bar{\mathbf{Z}}_n:=2\diag(\boldsymbol{\alpha}) \mathbf{Z}_n^* \diag(\boldsymbol{\alpha}^*)$; $\boldsymbol{\alpha}:=[1~\alpha~\alpha^2]^{\top}$; $\alpha=e^{-j\frac{2\pi}{3}}$; and $^*$ denotes complex conjugation. When not all phases are present, power injection and flow vectors and phase impedance matrices are zero-padded. For \eqref{eq:m-amodelv} to hold, the entries of $\mathbf{v}_n$ associated with non-existing phases are arbitrarily set to the corresponding entries of $\mathbf{v}_{\pi_n}$. 


\subsection{Network-Constrained EV Scheduling}\label{subsec:jo-model}
To facilitate network-constrained EV scheduling, the base active and reactive power loads $\{(\mathbf{d}_n(t),\mathbf{q}_n^d(t))\}$ for all $n$ and $t$ need to be predicted in advance. Active power loads $\mathbf{p}_n^d(t)$ consist of two parts: the base loads $\mathbf{d}_n(t)$ and the EV charging load. If $p_{n,\phi}^d(t)$ and $d_{n,\phi}(t)$ are respectively the total active load and the base load on phase $\phi$ of bus $n$, it holds that $p_{n,\phi}^d(t)=d_{n,\phi}(t)+\sum_{m\in\mathcal{M}_{n,\phi}}e_m(t)$. The cost $f_0(\mathbf{P}_0(t))$ of power supply from the main grid is convex and known in advance. Variables $\mathbf{p}_n^g(t)$ capture possible dispatchable generation distributed across the feeder, and $f_n^g(\mathbf{p}_n^g(t))$ is the associated convex quadratic cost for all $n\in\mathcal{N}$ and $t\in\mathcal{T}$.

To capture operational constraints, the following limits are introduced. Let $(\ubar{p}_{n,\phi}^g,\ubar{q}_{n,\phi}^g)$ be the lower, and $(\bar{p}_{n,\phi}^g,\bar{q}_{n,\phi}^g)$ the upper limits for distributed generation at phase $\phi\in\mathcal{P}_n$ of bus $n$. Define also $(\ubar{v}_{n,\phi},\bar{v}_{n,\phi})$ as the limits of squared voltage magnitudes at phase $\phi\in\mathcal{P}_n$ of bus $n$, $\bar{S}_n$ as the apparent power flow limits on line $n$, and $\bar{S}_f$ as the rated capacity of the feeder transformer. The utility company aims to minimize the total operation cost by coordinating vehicle charging and generation dispatch, while respecting charging and operational limitations. The pertinent network-constrained EV scheduling task can be posed as:
\allowdisplaybreaks[4] 
\begin{subequations}\label{eq:opf-ev}
\begin{align} 
\min ~&  \sum_{t \in \mathcal{T}}\left[f_0(\mathbf{P}_0(t)) + \sum_{n\in\mathcal{N}} f_n^g\left(\mathbf{p}_n^g(t)\right)\right] \label{eq:uobj}\\
\textrm{over}~&\{\mathbf{p}_n^g(t),\mathbf{q}_n^g(t),\mathbf{P}_n(t),\mathbf{Q}_n(t), \mathbf{v}_n(t)\}_{n\in\mathcal{N},t\in\mathcal{T}}, \{\mathbf{e}_m\}\nonumber\\
\textrm{s.to}~ 
~&\mathbf{p}_n^g(t)-\mathbf{p}_n^d(t) =\sum_{k\in\mathcal{C}_n}\mathbf{P}_k(t)  - \mathbf{P}_n(t),~\forall~n,t\label{eq:nodalpower}\\
~&\mathbf{q}_n^g(t)-\mathbf{q}_n^d(t)=\sum_{k\in\mathcal{C}_n}\mathbf{Q}_k(t)  - \mathbf{Q}_n(t),~\forall~n,t\label{eq:nodalrpower}\\
~&\mathbf{v}_{\pi_n}(t) - \mathbf{v}_n(t) =  \real\{\bar{\mathbf{Z}}_n(\mathbf{P}_n(t)+j\mathbf{Q}_n(t))\},\forall n,t\label{eq:voltage}\\
~&\ubar{p}_{n,\phi}^g\leq p_{n,\phi}^g(t) \leq \bar{p}_{n,\phi}^g,~\forall~\phi \in \mathcal{P}_n,n,t\label{eq:plimits}\\
~&\ubar{q}_{n,\phi}^g\leq q_{n,\phi}^g(t) \leq \bar{q}_{n,\phi}^g,~\forall~\phi \in \mathcal{P}_n,n,t\label{eq:qlimits}\\
~&\ubar{v}_n\leq v_{n,\phi}(t) \leq \bar{v}_n, ~\forall~\phi \in \mathcal{P}_n,n,t\label{eq:vlimits}\\
~& P_{n,\phi}^2(t)+Q_{n,\phi}^2(t) \leq \bar{S}_n^2,~\forall~\phi \in \mathcal{P}_n,n\in\mathcal{N}^+,t \label{eq:network}\\  
~& p_{n,\phi}^d(t)=  d_{n,\phi}(t) + \sum_{m\in\mathcal{M}_{n,\phi}}e_m(t), \forall~ \phi\in\mathcal{P}_n,n,t \label{eq:coupledload}\\
~&\mathbf{e}_m \in\mathcal{E}_m,~\forall~m\label{eq:ev}\\
~&(\mathbf{1}^\top\mathbf{P}_0(t))^2+(\mathbf{1}^\top\mathbf{Q}_0(t))^2 \leq \bar{S}_f^2,~\forall~ t.\label{eq:feeder0}
\end{align}
\end{subequations}
Constraints \eqref{eq:nodalpower}--\eqref{eq:voltage} originate from the power flow model; constraints \eqref{eq:plimits}--\eqref{eq:qlimits} enforce generation limits; voltage regulation is guaranteed via \eqref{eq:vlimits}; apparent power flows are upper bounded by \eqref{eq:network}; the equalities in \eqref{eq:coupledload} define demands across phases and buses; constraint \eqref{eq:ev} is related to the per-vehicle charging profile; and \eqref{eq:feeder0} results from the capacity limit of the feeder transformer. 

The cost functions and all the constraints apart from the EV charging constraint in \eqref{eq:ev} are \emph{separable across time}. The capacity limit in \eqref{eq:feeder0} couples flows across phases, while the voltage regulation constraints in \eqref{eq:voltage} and \eqref{eq:vlimits}  couple variables across buses and phases. For linear and convex quadratic costs, problem \eqref{eq:opf-ev} can be reformulated as a standard quadratically-constrained quadratic program and tackled by standard solvers in a centralized manner. Nonetheless, for increasing grid sizes, longer time horizons $\mathcal{T}$, and/or shorter control periods, tackling \eqref{eq:opf-ev} could be challenging. In addition, private information on a per-vehicle basis needs to be collected and processed by the utility. These considerations motivate well the privacy-preserving and scalable (both in space and time) scheme for solving \eqref{eq:opf-ev} that is pursued next. 

\section{Distributed Optimal Charging Protocol}\label{sec:jo-admm}
This section delineates an ADMM-based method for decomposing \eqref{eq:opf-ev} into smaller subproblems. Notably, each subproblem either enjoys a closed-form solution or it can be tackled efficiently by Alg.~\ref{alg:docsev}. As a brief review, ADMM solves problems of the form~\cite{Boyd10}
\begin{align}\label{eq:ADMM}
\min_{\mathbf{x}\in\mathcal{X},\mathbf{z}\in\mathcal{Z}} & \left\{f(\mathbf{x}) + g(\mathbf{z}):~\mathbf{F}\mathbf{x} + \mathbf{G}\mathbf{z}=\mathbf{b}\right\}
\end{align}
where $f(\mathbf{x})$ and $g(\mathbf{z})$ are convex functions; $\mathcal{X}$ and $\mathcal{Z}$ are convex sets; and $(\mathbf{F},\mathbf{G},\mathbf{b})$ model the linear equality constraints coupling variables $\mathbf{x}$ and $\mathbf{z}$. In its normalized form, ADMM assigns a Lagrange multiplier $\mathbf{w}$ for the equality constraint and solves \eqref{eq:ADMM} by iterating over the following three recursions
\begin{subequations}\label{eq:ADMM:steps}
\begin{align}
\mathbf{x}^{i+1}&\in\arg\min_{\mathbf{x}\in\mathcal{X}} f(\mathbf{x}) + \tfrac{\rho}{2}\|\mathbf{F}\mathbf{x} + \mathbf{G}\mathbf{z}^i-\mathbf{b} + \mathbf{w}^i\|_2^2\label{eq:ADMM:S1}\\
\mathbf{z}^{i+1}&\in\arg\min_{\mathbf{z}\in\mathcal{Z}} g(\mathbf{z}) + \tfrac{\rho}{2}\|\mathbf{F}\mathbf{x}^{i+1} + \mathbf{G}\mathbf{z}-\mathbf{b} + \mathbf{w}^i\|_2^2\label{eq:ADMM:S2}\\
\mathbf{w}^{i+1}&:=\mathbf{w}^{i} + \mathbf{F}\mathbf{x}^{i+1} + \mathbf{G}\mathbf{z}^{i+1}-\mathbf{b}\label{eq:ADMM:S3}
\end{align}
\end{subequations}
for some $\rho>0$. ADMM has been successfully applied to decentralize various power system tasks across buses~\cite{HaozhuD},~\cite{KeGi12},~\cite{QiuyuPeng}. Related ideas are adopted here to decouple the spatially-coupled constraints \eqref{eq:nodalpower}--\eqref{eq:voltage}. 

To that end, each bus $n\in\mathcal{N}$ maintains a local copy of the variables associated with the squared voltage magnitude of its parent bus, and the power flows feeding its children buses. These auxiliary variables are marked with a hat as $\hat{\mathbf{v}}_{n}$ and $\{(\hat{\mathbf{P}}_{k},\hat{\mathbf{Q}}_{k})\}_{k\in\mathcal{C}_n}$. The duplicate variable $\hat{\mathbf{v}}_{n}$ stored at bus $n$ should agree with the original variable $\mathbf{v}_{\pi_n}$ stored at bus $\pi_n$. To decentralize the computations, we further introduce the \emph{consensus} variable $\tilde{\mathbf{v}}_{\pi_n}$, and impose the constraints $\mathbf{v}_{\pi_n}=\tilde{\mathbf{v}}_{\pi_n}$ and $\hat{\mathbf{v}}_{n}=\tilde{\mathbf{v}}_{\pi_n}$ for all non-leaf buses. By repeating this process for power flow variables and all $n\in\mathcal{N}$, the physical grid model will be later decoupled across buses.

We also introduce duplicate variables $\{\tilde{\mathbf{p}}_n^d(t)\}_{n\in\mathcal{N}}$ for net loads to separate the tasks of EV charging and generation dispatch. As detailed later, imposing the constraints $\tilde{\mathbf{p}}_n^d(t)=\mathbf{p}_n^d(t)$ for all $n$, enables isolating \eqref{eq:ev} from the rest of the constraints in \eqref{eq:opf-ev}; resulting in localized EV charging subproblems that is a special case of~\eqref{eq:EVS2}. 


For a compact representation define the aggregate variables:
\begin{align*}
\mathbf{x}_n(t)&:=\left\{\mathbf{v}_{n}(t),\mathbf{p}_n^{g}(t),\mathbf{p}_n^{d}(t),\mathbf{q}_{n}^{g}(t),\mathbf{P}_{n}(t),\mathbf{Q}_{n}(t)\right\}\\
\hat{\mathbf{x}}_n(t)&:=\left\{\hat{\mathbf{v}}_{n}(t),\{\hat{\mathbf{P}}_{k}(t),\hat{\mathbf{Q}}_{k}(t)\}_{k\in \mathcal{C}_n}\right\}\\
\tilde{\mathbf{z}}_n(t)&:=\left\{\tilde{\mathbf{v}}_{n}(t),\tilde{\mathbf{p}}_{n}^{g}(t),\tilde{\mathbf{p}}_{n}^{d}(t),\tilde{\mathbf{q}}_{n}^{g}(t),\tilde{\mathbf{P}}_{n}(t),\tilde{\mathbf{Q}}_{n}(t)\right\}
\end{align*}
for all $n\in\mathcal{N}$ and $t\in\mathcal{T}$. With the newly introduced variables, problem \eqref{eq:opf-ev} can be equivalently expressed as:
\begin{subequations}\label{eq:augadmm}
\begin{align} \label{eq:augutility}
\min ~&  \sum_{t \in \mathcal{T}}\left[f_0(\mathbf{P}_0(t)) + \sum_{n\in\mathcal{N}} f_n^g(\tilde{\mathbf{p}}_n^g(t))\right]\\
\textrm{over}~&\{\mathbf{x}_n(t),\hat{\mathbf{x}}_n(t),\tilde{\mathbf{z}}_n(t)\}_{n\in\mathcal{N},t\in\mathcal{T}},~\{\mathbf{e}_m\in\mathcal{E}_m\}_{m\in\mathcal{M}}, \nonumber\\
\textrm{s.to}~& \mathbf{p}_n^g(t)-\mathbf{p}_n^d(t) =\sum_{k\in\mathcal{C}_n}\hat{\mathbf{P}}_k(t)  - \mathbf{P}_n(t),\forall~ n\in\mathcal{N},t \label{eq:augp} \\
~&\mathbf{q}_n^g(t)-\mathbf{q}_n^d(t)=\sum_{k\in\mathcal{C}_n}\hat{\mathbf{Q}}_k(t)  - \mathbf{Q}_n(t),\forall~  n\in\mathcal{N},t \label{eq:augq}\\
~&\hat{\mathbf{v}}_{n}(t) - \mathbf{v}_n(t) = \real\{\bar{\mathbf{Z}}_n(\mathbf{P}_n(t)+j\mathbf{Q}_n(t))\}\nonumber\\
~&\quad \quad\forall n\in\mathcal{N}^+,t \label{eq:augv}\\
~&\ubar{p}_{n,\phi}^g\leq \tilde{p}_{n,\phi}^g(t) \leq \bar{p}_{n,\phi}^g,~\forall~ \phi \in \mathcal{P}_n,n\in\mathcal{N},t \label{eq:augpg} \\
~&\ubar{q}_{n,\phi}^g\leq \tilde{q}_{n,\phi}^g(t) \leq \bar{q}_{n,\phi}^g,~\forall~ \phi \in \mathcal{P}_n,n\in\mathcal{N},t\\
~&\ubar{v}_n\leq \tilde{v}_{n,\phi}(t) \leq \bar{v}_n,~\forall~ \phi \in \mathcal{P}_n,n\in\mathcal{N},t\\
~&\tilde{P}_{n,\phi}^2(t)+\tilde{Q}_{n,\phi}^2(t) \leq \bar{S}_n^2,~\forall~ \phi \in \mathcal{P}_n,n\in \mathcal{N}_+,t \label{eq:augS}\\  
~&\mathbf{P}_{n}(t) = \tilde{\mathbf{P}}_n(t),~\mathbf{Q}_{n}(t) = \tilde{\mathbf{Q}}_n(t),~\mathbf{v}_{n}(t)=\tilde{\mathbf{v}}_{n}(t), \nonumber\\
~&\quad\quad\forall~ n\in \mathcal{N}^+,t \label{eq:augadmm1}\\
~&\hat{\mathbf{P}}_{n}(t) = \tilde{\mathbf{P}}_n(t),~\hat{\mathbf{Q}}_{n}(t) = \tilde{\mathbf{Q}}_n(t),~\hat{\mathbf{v}}_{n}(t) = \tilde{\mathbf{v}}_{\pi_n}(t),\nonumber\\
~&\quad\quad\forall~ n\in \mathcal{N}^+,t \label{eq:augadmm2}\\
~&\mathbf{p}_{n}^g(t) = \tilde{\mathbf{p}}_n^g(t),\mathbf{p}_{n}^d(t)=\tilde{\mathbf{p}}_{n}^d(t),\mathbf{q}_{n}^g(t) = \tilde{\mathbf{q}}_n^g(t),\nonumber\\
~&\quad\quad\forall~ n\in\mathcal{N},t \label{eq:augadmm3}\\
~&\tilde{p}_{n,\phi}^{d}(t)=  d_{n,\phi}(t)+\sum_{m\in \mathcal{M}_{n,\phi}} e_m(t),~\forall \phi \in \mathcal{P}_{n}, n,t  \label{eq:augev}\\
~& (\mathbf{1}^\top \tilde{\mathbf{P}}_0(t))^2+(\mathbf{1}^\top\tilde{\mathbf{Q}}_0(t))^2 \leq \bar{S}_f^2,~\forall t  \label{eq:augf}
\end{align}
\end{subequations}
The equality constraints between duplicate variables in \eqref{eq:augadmm1}--\eqref{eq:augev} are assigned Langrange multipliers according to Table~\ref{table:lmultiplier}. Adopting the ADMM iterates of \eqref{eq:ADMM:steps} to solve \eqref{eq:augadmm}, variables $\{\mathbf{x}_n(t), \hat{\mathbf{x}}_n(t)\}_{n\in\mathcal{N},t\in\mathcal{T}}$ and $\{\mathbf{e}_m\}_{m\in\mathcal{M}}$ are updated in the first ADMM step, whereas variables $\left\{\{\tilde{\mathbf{z}}_n(t)\}_{n\in\mathcal{N}}\right\}_{t\in\mathcal{T}}$ are updated during the second ADMM step as detailed next.


\begin{table}[t]
\renewcommand{\arraystretch}{1.2}
\vspace*{-0.5em}
\caption{Lagrange multipliers for problem \eqref{eq:augadmm}}
\vspace*{-0.5em}
\centering 
\begin{tabular}{|c|c||c|c|} 
\hline
$\mathbf{p}_n^{g}(t)= \tilde{\mathbf{p}}_n^{g}(t)$ & $\boldsymbol{\lambda}_n^p(t)$& $\mathbf{q}_n^{g}(t)=\tilde{\mathbf{q}}_n^{g}(t)$ & $\boldsymbol{\lambda}_n^q(t)$\\
\hline 
$\hat{\mathbf{P}}_{n}(t) = \tilde{\mathbf{P}}_n(t)$& $\hat{\boldsymbol{\lambda}}_n^P(t)$ & $\mathbf{P}_{n}(t) = \tilde{\mathbf{P}}_n(t)$& $\boldsymbol{\lambda}_n^P(t)$ \\ 
\hline
$\hat{\mathbf{Q}}_{n}(t) = \tilde{\mathbf{Q}}_n(t)$& $\hat{\boldsymbol{\lambda}}_n^Q(t)$ & $\mathbf{Q}_{n}(t) = \tilde{\mathbf{Q}}_n(t)$& $\boldsymbol{\lambda}_n^Q(t)$ \\
\hline
$\hat{\mathbf{v}}_{n}(t) = \tilde{\mathbf{v}}_{\pi_n}(t)$& $\hat{\boldsymbol{\lambda}}_n^v(t)$ & $\mathbf{v}_{n}(t) = \tilde{\mathbf{v}}_n(t)$& $\boldsymbol{\lambda}_n^v(t)$ \\
\hline
$\mathbf{p}_n^{d}(t)= \tilde{\mathbf{p}}_n^{d}(t)$ & $\boldsymbol{\lambda}_n^d(t)$& Constraints \eqref{eq:augev}& $\mu_{n,\phi}(t)$ \\
\hline 
\end{tabular}
\label{table:lmultiplier}
\vspace*{-1.5em}
\end{table}

\subsection{First Step of ADMM}\label{subsec:S1}
Due to the form the generic update \eqref{eq:ADMM:S1} takes for the problem at hand, variables $\{\mathbf{x}_n(t),\hat{\mathbf{x}}_n(t)\}_{n\in\mathcal{N},t\in\mathcal{T}}$ can be updated separately from the EV charging profiles $\{\mathbf{e}_m\}_{m\in\mathcal{M}}$. The updates for these two variable sets are studied next.

Heed that $\{\mathbf{x}_n(t),  \hat{\mathbf{x}}_n(t)\}_{n\in\mathcal{N},t\in\mathcal{T}}$ can be optimized independently across buses and time periods. Nevertheless, for fixed bus and time indices $(n,t)$, variables $\mathbf{x}_n(t)$ and $\hat{\mathbf{x}}_n(t)$ are coupled due to constraints \eqref{eq:augp}--\eqref{eq:augv}. To simplify the presentation, we drop the time index and consider the canonical subproblems involved for all $t\in\mathcal{T}$. Let $\hat{\mathbf{z}}_n:= \left\{\tilde{\mathbf{v}}_{\pi_n},\{\tilde{\mathbf{P}}_{k},\tilde{\mathbf{Q}}_{k}\}_{k\in \mathcal{C}_n}\right\}$ for bus $n\in\mathcal{N}^+$. Variables $\mathbf{x}_n$ and $\hat{\mathbf{x}}_n$ are updated during the $i$-th iteration as the minimizers of
\begin{align} \label{eq:xnodeopf}
\min_{\mathbf{x}_n,\hat{\mathbf{x}}_n} ~&~\|\mathbf{x}_n-\mathbf{z}_n^i+\boldsymbol{\lambda}_n^i\|_2^2+\|\hat{\mathbf{x}}_n -\hat{\mathbf{z}}_n^i+\hat{\boldsymbol{\lambda}}_n^i\|_2^2\nonumber\\ \textrm{s.to}~&~\eqref{eq:augp}-\eqref{eq:augv}. 
\end{align}
For $n=0$ and due to the power supply cost from the main grid, variables $(\mathbf{P}_0,\mathbf{Q}_0)$ are found as the minimizers of 
\begin{align}\label{eq:xsubstation}
\min_{\mathbf{P}_0,\mathbf{Q}_0}~&~\|\mathbf{x}_0-\mathbf{z}_0^i+\boldsymbol{\lambda}_0^i\|_2^2+\|\hat{\mathbf{x}}_0 -\hat{\mathbf{z}}_0^i+\hat{\boldsymbol{\lambda}}_0^i\|_2^2 +\frac{2}{\rho}f_0({\mathbf{P}}_0)\nonumber \\
\textrm{s.to}~ &~\eqref{eq:augp}-\eqref{eq:augq}. 
\end{align} 
Problems \eqref{eq:xnodeopf}--\eqref{eq:xsubstation} are linearly-constrained quadratic programs with closed-form minimizers~\cite{BoVa04}. 

We next focus on updating the vehicle charging profiles $\{\mathbf{e}_m\}_{m\in\mathcal{M}}$ at iteration $i$. Interestingly, the task of EV charging decouples over buses and phases. The charging profiles for vehicles $m\in\mathcal{M}_{n,\phi}$ can be updated as the minimizers of
\begin{equation}\label{eq:EVS}
\min_{\{\mathbf{e}_m \in\mathcal{E}_m\}_{m\in\mathcal{M}_{n,\phi}}} \frac{1}{2}\sum_{t\in \mathcal{T}} \left(l_{n,\phi}^i(t)+\sum_{m \in \mathcal{M}_{n,\phi}}e_m(t)\right)^2
\end{equation}
where $l_{n,\phi}^i(t):=d_{n,\phi}(t)-\tilde{p}_{n,\phi}^{d,i}(t)+{\mu}_{n,\phi}^i(t)$ and the Lagrange multiplier ${\mu}_{n,\phi}^i(t)$  reflects the network constraints. Note \eqref{eq:EVS} is actually a special case of~\eqref{eq:EVS2} with $C_t(x)=x^2/2,~ \forall t$. Hence, subproblem \eqref{eq:EVS} can be solved using Alg.~\ref{alg:docsev}. 

\begin{figure}[t]
    \centering
    \subfloat[First step of ADMM.]{\includegraphics[width=0.243\textwidth]{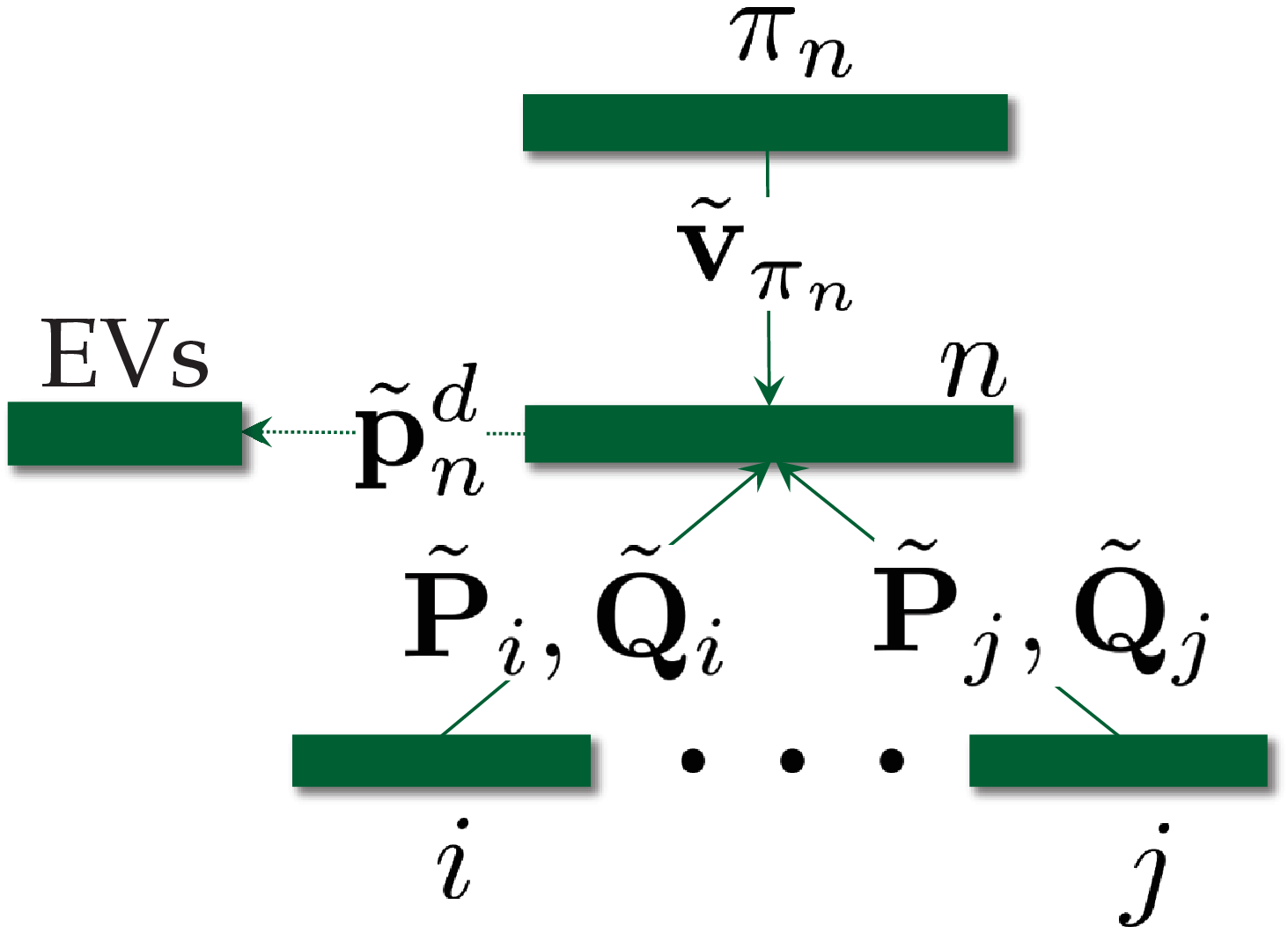}\label{fig:admm1}} 
    ~ 
   \subfloat[Second step of ADMM. ]{    \includegraphics[width=0.241\textwidth]{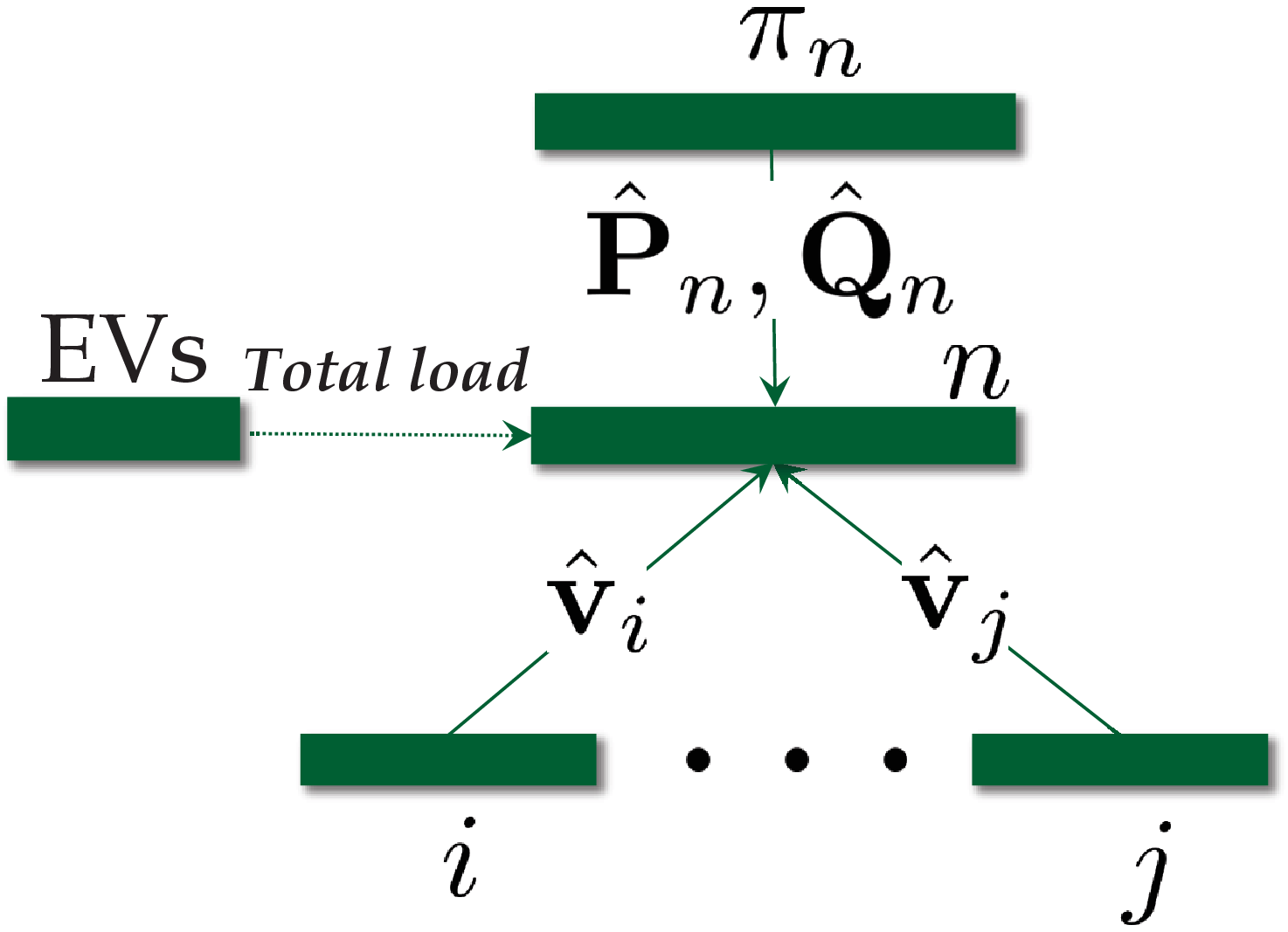}\label{fig:admm2}}
    \caption{Information exchange in the ADMM steps for bus $n$.} \label{fig:comm}
\end{figure}

In the first step of ADMM, each bus $n$ needs to collect $\tilde{\mathbf{v}}_{\pi_n}$ from its parent and $\{\tilde{\mathbf{P}}_{k},\tilde{\mathbf{Q}}_{k}\}_{k\in \mathcal{C}_n}$ from all its children as depicted in Fig.~\ref{fig:admm1}.  Meanwhile, each bus $n$ transfers $\tilde{\mathbf{p}}_{n}^d$ to its EV scheduling center, where the charging profile of EVs are optimized using Alg.~\ref{alg:docsev}.    
%

\subsection{Second Step of ADMM}\label{subsec:S2}
Finding optimal $\{\tilde{\mathbf{z}}_n(t)\}_{n\in\mathcal{N},t\in\mathcal{T}}$ can be performed independently across buses and time slots. Because of that, time indices are ignored. Every bus $n$ has to solve five sub-problems in parallel, each one associated with the variables $\tilde{\mathbf{p}}_{n}^{g}$,  $\tilde{\mathbf{p}}_{n}^{d}$, $\tilde{\mathbf{q}}_{n}^{g}$, $\tilde{\mathbf{v}}_{n}$, and $(\tilde{\mathbf{P}}_{n},\tilde{\mathbf{Q}}_{n})$. Firstly, updating $\tilde{\mathbf{p}}_{n}^{g}$, $\tilde{\mathbf{p}}_{n}^{d}$, $\tilde{\mathbf{q}}_{n}^{g}$, and $\tilde{\mathbf{v}}_{n}$ decouples over phases of bus $n$ too. It can be shown that per phase variables are updated as the minimizers of a univariate convex quadratic function possibly over box constraints. If the generation cost at bus $n$ is $f_n^g(\tilde{\mathbf{p}}_n^g):=\sum_{\phi \in \mathcal{P}_n} a_{n,\phi}({p}_{n,\phi}^{g})^2+b_{n,\phi}{p}_{n,\phi}^{g}+c_{n,\phi}$ with $a_{n,\phi}\geq 0$, then $\tilde{p}_{n,\phi}^{g}$ is updated at iteration $i$ by solving
\begin{align} \label{eq:znodeopfq}
\min_{\tilde{p}_{n,\phi}^{g}} ~&~ a_{n,\phi}(\tilde{p}_{n,\phi}^{g})^2+b_{n,\phi}\tilde{p}_{n,\phi}^{g} + \tfrac{\rho}{2}(p_{n,\phi}^{g,i}- \tilde{p}_{n,\phi}^{g}+\lambda_{n,\phi}^{p,i})^2\nonumber\\
\textrm{s.to}~&~\ubar{p}_{n,\phi}^g\leq \tilde{p}_{n,\phi}^g \leq \bar{p}_{n,\phi}^g. 
\end{align}
The minimizer of \eqref{eq:znodeopfq} is expressed as
\begin{equation}\label{eq:oppg}
\tilde{p}_{n,\phi}^{g,i+1}=\left[\frac{\rho(p_{n,\phi}^{g,i}+\lambda_{n,\phi}^{p,i})-b_{n,\phi}}{2a_{n,\phi}+\rho}\right]_{\ubar{p}_{n,\phi}^g}^{\bar{p}_{n,\phi}^g}
\end{equation}
where $[x]_{\ubar{x}}^{\bar{x}}:=\max\{\ubar{x},\min\{x,\bar{x}\}\}$. The entries of $\tilde{\mathbf{q}}_n^{g}$ and $\tilde{\mathbf{v}}_n$ are similarly found as
\begin{align}
\tilde{q}_{n,\phi}^{g,i+1}&=\left[{q}_{n,\phi}^{g,i}+\lambda_{n,\phi}^{q,i}\right]_{\ubar{q}_{n,\phi}^g}^{\bar{q}_{n,\phi}^g}\label{eq:opqg}\\
\tilde{v}_{n,\phi}^{i+1}&=\left[\frac{\sum_{k\in \mathcal{C}_n}(\hat{{v}}_{k,\phi}^i+\hat{\lambda}_{k,\phi}^{v,i})+v_{n,\phi}^i+\lambda_{n,\phi}^{v,i}}{|\mathcal{C}_n|+1} \right]_{\ubar{v}_{n}}^{\bar{v}_{n}}.\label{eq:opv}
\end{align}
The entries of $\tilde{\mathbf{p}}_n^{d}$ are obtained as the solutions of unconstrained univariate convex quadratic programs as
\begin{equation}
\tilde{p}_{n,\phi}^{d,i+1}=\frac{1}{2} \left({p}_{n,\phi}^{d,i} + \lambda_{n,\phi}^{d,i} + d_{n,\phi}
+\sum_{m \in \mathcal{M}_{n,\phi}} e_m^i + \mu_{n,\phi}^i\right).\label{eq:oppd}
\end{equation} 

The optimizations involved in updating the consensus power flow variables $\{\tilde{\mathbf{P}}_n,\tilde{\mathbf{Q}}_n\}_{n\in\mathcal{N}^+}$ decouple across phases. The consensus power flow variables $\{(\tilde{P}_{n,\phi},\tilde{Q}_{n,\phi})\}_{\phi\in\mathcal{P}_n,n\in\mathcal{N}^+}$ are updated by solving the problems for all $\phi\in\mathcal{P}_n$ and $n\in\mathcal{N}^+$:
\begin{align} \label{eq:PQ}
\min_{\tilde{P}_{n,\phi},\tilde{Q}_{n,\phi}} ~&(\tilde{P}_{n,\phi}-\breve{P}_{n,\phi}^i)^2 + (\tilde{Q}_{n,\phi}-\breve{Q}_{n,\phi}^i)^2\\
\textrm{s.to}~~&\tilde{P}_{n,\phi}^2+\tilde{Q}_{n,\phi}^2 \leq \bar{S}_n^2\nonumber
\end{align}
where $\breve{P}_{n,\phi}^i:=\tfrac{1}{2}(P_{n,\phi}^i+\hat{P}_{n,\phi}^i + \lambda_{n,\phi}^{P,i} + \hat{\lambda}_{n,\phi}^{P,i})$, and $\breve{Q}_{n,\phi}^i:=\tfrac{1}{2}(Q_{n,\phi}^i+\hat{Q}_{n,\phi}^i + \lambda_{n,\phi}^{Q,i} + \hat{\lambda}_{n,\phi}^{Q,i})$. Resorting to the KKT conditions for \eqref{eq:PQ} shows that its minimizers are
\begin{subequations}\label{eq:PQclosedform}
\begin{align}
\tilde{P}_{n,\phi}^{i+1}&:=\min\left\{\frac{\bar{S}_n}{\sqrt{(\breve{P}_{n,\phi}^i)^2+(\breve{Q}_{n,\phi}^i)^2}},1\right\}\breve{P}_{n,\phi}^i\label{eq:Pclosedform}\\
\tilde{Q}_{n,\phi}^{i+1}&:=\min\left\{\frac{\bar{S}_n}{\sqrt{(\breve{P}_{n,\phi}^i)^2+(\breve{Q}_{n,\phi}^i)^2}},1\right\}\breve{Q}_{n,\phi}^i.\label{eq:Qclosedform}
\end{align}
\end{subequations}
The substation power flows are updated as the solution to
\begin{align} \label{eq:PQ0}
\min_{\tilde{\mathbf{P}}_0,\tilde{\mathbf{Q}}_0} ~&\|\tilde{\mathbf{P}}_0-\breve{\mathbf{P}}_0^i\|_2^2 +\|\tilde{\mathbf{Q}}_0-\breve{\mathbf{Q}}_0^i\|_2^2\\ 
\textrm{s.to}~~& (\mathbf{1}^{\top}\tilde{\mathbf{P}}_0)^2 + (\mathbf{1}^{\top}\tilde{\mathbf{Q}}_0)^2 \leq  \bar{S}_f^2\nonumber
\end{align}
where $\breve{\mathbf{P}}_0^i:=\mathbf{P}_0^i+\boldsymbol{\lambda}_{0}^{P,i}$ and $\breve{\mathbf{Q}}_0^i:=\mathbf{Q}_0^i+\boldsymbol{\lambda}_{0}^{Q,i}$. The following optimal solution to~\eqref{eq:PQ0} is derived in the Appendix 
\begin{proposition} \label{pro:2}
The optimal solution of problem \eqref{eq:PQ0} is
\begin{subequations}\label{eq:prop}
\begin{align}
\tilde{\mathbf{P}}_0^{i+1}&:=\breve{\mathbf{P}}_0^i-\max\left\{1-\frac{\bar{S}_f}{\Sigma},0 \right\} \frac{\mathbf{1}\mathbf{1}^{\top} \breve{\mathbf{P}}_0^i}{3} \label{eq:pro2P}\\
\tilde{\mathbf{Q}}_0^{i+1}&:=\breve{\mathbf{Q}}_0^i-\max\left\{1-\frac{\bar{S}_f}{\Sigma},0 \right\}\frac{\mathbf{1}\mathbf{1}^{\top} \breve{\mathbf{Q}}_0^i}{3}\label{eq:pro2Q}
\end{align}
\end{subequations}
where $\Sigma:=\sqrt{(\mathbf{1}^{\top} \breve{\mathbf{P}}_0^i)^2+(\mathbf{1}^{\top} \breve{\mathbf{Q}}_0^i)^2}$.
\end{proposition}

To implement the second step of ADMM,  bus $n$ gathers its copies $(\hat{\mathbf{P}}_{n},\hat{\mathbf{Q}}_{n})$ from its parent, $\{\hat{v}_k\}_{k\in \mathcal{C}_n}$ from all its children, and the total charging load $\{\sum_{m\in \mathcal{M}_{n,\phi}} \mathbf{e}_m^{k}\}_{\phi \in \mathcal{P}_n}$ of all the connected EVs as presented in Fig.~\ref{fig:admm2}. Then bus $n$ updates $\tilde{\mathbf{z}}_n$ according to \eqref{eq:znodeopfq}--\eqref{eq:oppd}, \eqref{eq:PQclosedform}, and \eqref{eq:prop}. 

The Lagrange multipliers are updated according to~\eqref{eq:ADMM:S3}, i.e., every multiplier is equal to its previous value plus the most recent constraint violation. 

\section{Numerical Tests}\label{sec:tests}

\subsection{Frank-Wolfe Scheme for Vehicle Charging} \label{subsec:simua1}
We first evaluated Alg.~\ref{alg:docsev} by simulating the charging of 59  EVs. The costs are selected as $C_t(x):=x^2/2$ for all $t$~\cite{GTL11}. For all vehicles, the battery capacity was 20kWh and the maximum charging rate was 3.45kW~\cite{NCEV}. The plug-in/-out times and daily travel miles were set according to the statistical estimates from survey travel data~\cite{Federalhighway}, \cite{EVload}.
The expected state of charge for EVs was fixed to 90\%, and the energy needed per 100km is $E_{100}$=15kWh. The initial state of charge for EV $m$ was modeled as $s_m^{soc}=0.9-D^{\text{miles}}_m E_{100}/(100B_m)$ for daily travel miles $D^{\text{miles}}_m$. Normalized base load curves were obtained by averaging the 2014 residential load data from Southern California Edison. A day-long horizon starting at midnight was divided into $T=96$ slots. Tests were run on Matlab using an Intel CPU @ 3.6 GHz (32 GB RAM) computer.

\begin{figure}[t]
\centering
\includegraphics[scale=0.3]{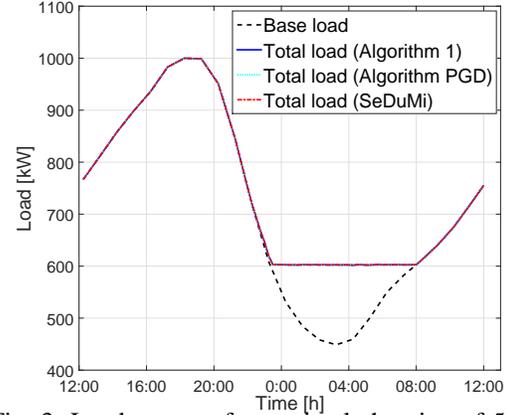}
\vspace*{-1em}
\caption{Load curves after optimal charging of 59 EVs.}
\vspace*{-1em}
\label{fig:fwodc}
\end{figure}

Parameter $d(t)$ was the normalized residential load with the maximum load set to 1000 kW~\cite{GTL11}. The minimizer of \eqref{eq:opf-ev} was obtained via SeDuMi, Algorithm~\ref{alg:docsev}, and the PGD solver of \cite{GTL11}. The subproblem~\eqref{eq:PGD} entailed in PGD was solved by SeDuMi. Algorithm~\ref{alg:docsev} and PGD were terminated once the relative cost error denoted as $\epsilon$ became smaller than $10^{-7}$. Figure~\ref{fig:fwodc} shows that the three resultant load curves coincide and feature a flat load valley. Performing the updates for Alg.~\ref{alg:docsev} and PGD sequentially in one computer, Alg.~\ref{alg:docsev} converged within 0.78 sec, PGD in 734.74 sec, and SeDuMi (centralized solver) in 82.47 sec. Had Alg.~\ref{alg:docsev} and PGD run in parallel, shorter running times would have been obtained. 

\begin{figure}[t]
\centering
\includegraphics[scale=0.3]{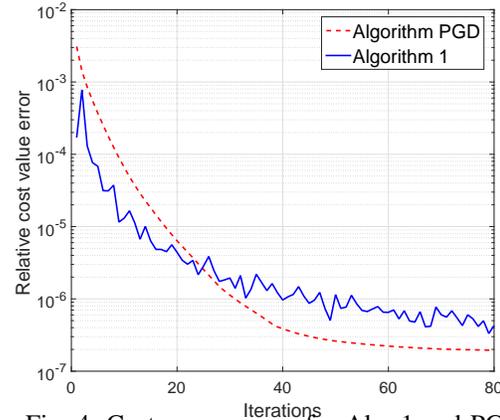}
\vspace*{-1em}
\caption{Cost convergence for Alg.~\ref{alg:docsev} and PGD.}
\vspace*{-1em}
\label{fig:con}
\end{figure}

\begin{figure}[t]
\centering
\includegraphics[scale=0.3]{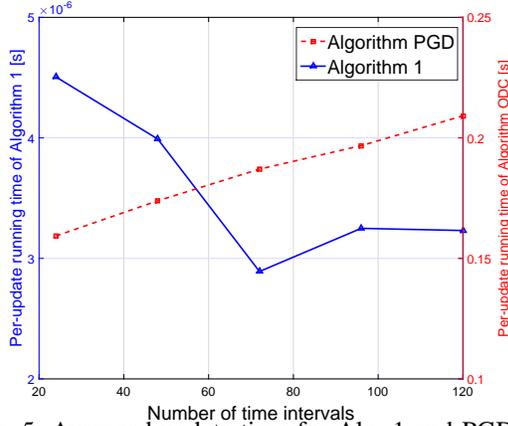}
\vspace*{-1em}
\caption{Averaged update time for Alg.~\ref{alg:docsev} and PGD.}
\vspace*{-1em}
\label{fig:runt}
\end{figure}

Figure~\ref{fig:con} depicts the cost convergence curves for Alg.~\ref{alg:docsev} and PGD for scheduling 34 vehicles. Setting $\epsilon=10^{-6}$, Alg.~\ref{alg:docsev} was terminated after 80 iterations, and PGD was also run for 80 iterations. Observed from Fig.~\ref{fig:con}, the decreasing rate of Alg.~\ref{alg:docsev} is almost the same as PGD; though the average update time for Alg.~\ref{alg:docsev} is $3.5 \times 10^{-6}$~sec, which is significantly superior to PGD's average update of 1.8~sec. Figure~\ref{fig:runt} presents the running time (averaged over all EVs and iterations) for a single update. It is worth stressing that Alg.~\ref{alg:docsev} requires roughly $10^{-6}$ sec regardless of the number of time intervals while PGD's average update time increases almost linearly with the number of time intervals. The major computational advantage of Alg.~\ref{alg:docsev} is the simple update in~\eqref{eq:opts}.

\subsection{ADMM-based Scheme for Network-Constrained Charging}\label{subsec:simua1ncev}
\begin{figure}[t]
    \centering
    \subfloat[Objective value convergence.]{\includegraphics[width=0.243\textwidth]{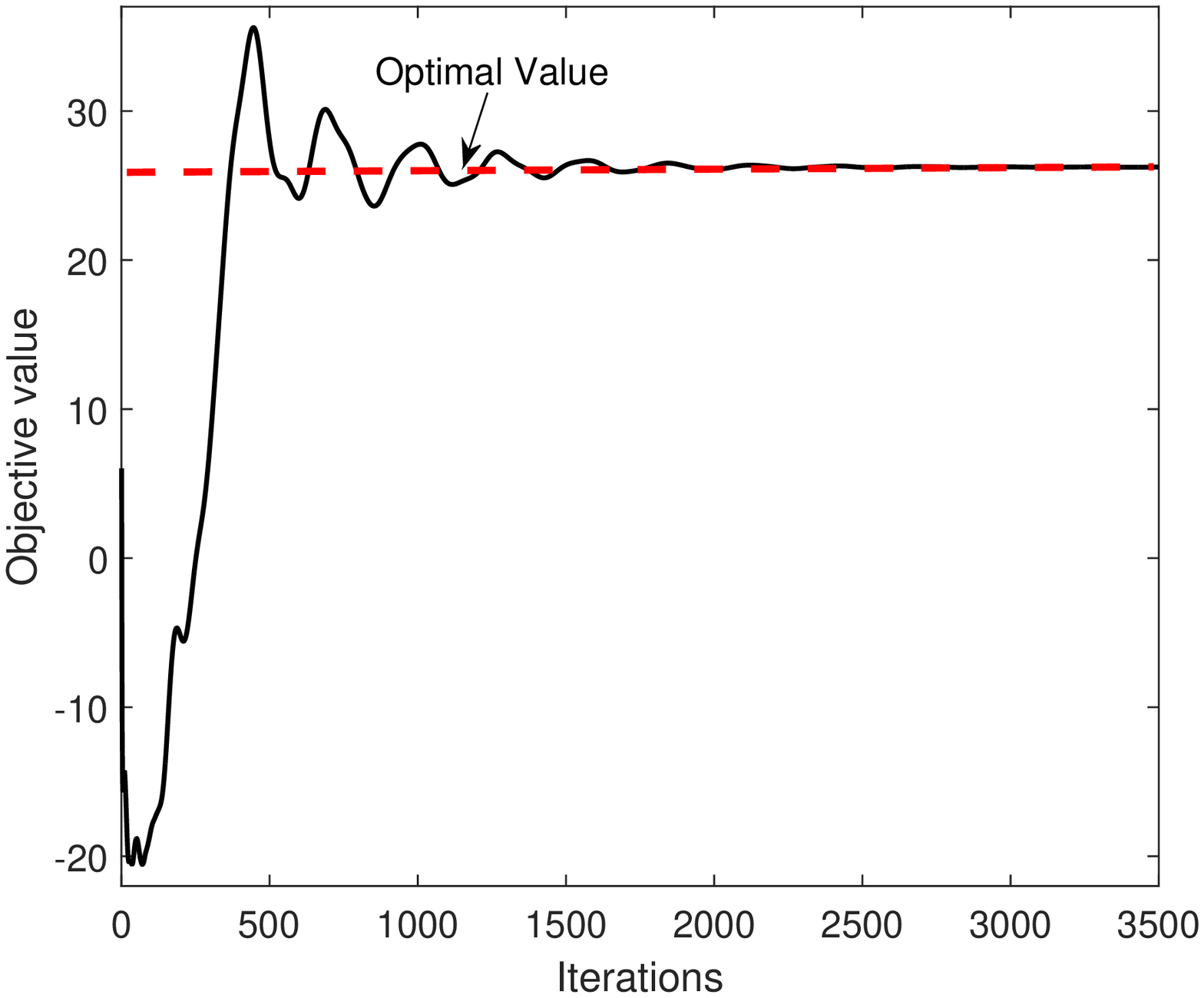}\label{fig:dp1}} 
    ~ 
   \subfloat[Primal/dual residual convergence.]{\includegraphics[width=0.241\textwidth]{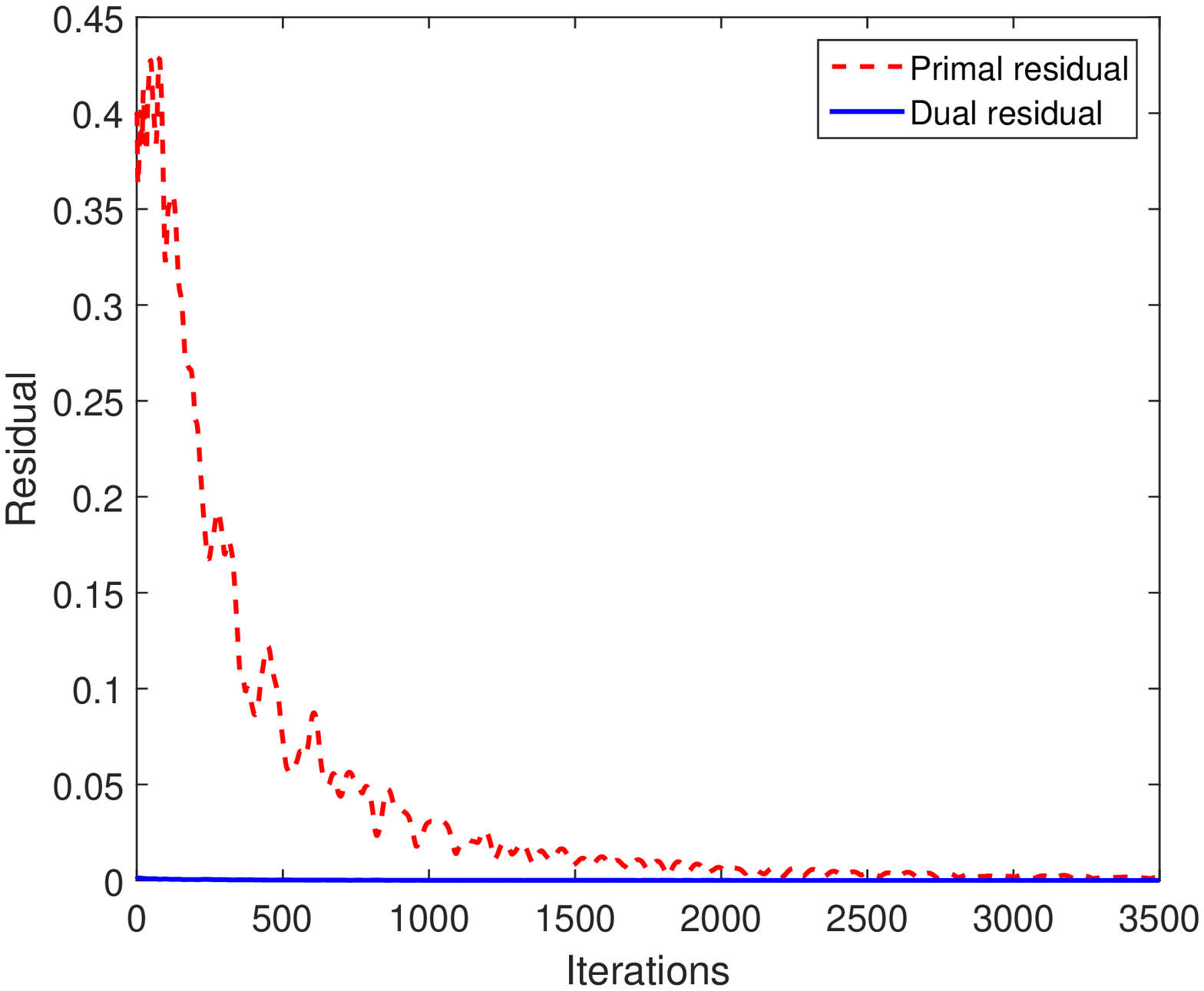}\label{fig:dp2}}
    \caption{Decentralized network-constrained EV charging.} \label{fig:dpncevs}
    \vspace*{-0.7em}
\end{figure}

The decentralized algorithm for network-constrained EV charging was tested using the unbalanced IEEE 123-bus feeder~\cite{TestFeeder}. Fifteen distributed generation (DG) units were located in the system; while 5, 10, 15, 25 and 5 EVs were being charged on bus 3, 15, 64, 82, and 102, respectively. At iteration $i$, the primal and dual residual for~\eqref{eq:ADMM:steps} are defined as $o_p^i:=\|\mathbf{F}\mathbf{x}^{i} + \mathbf{G}\mathbf{z}^{i}-\mathbf{b} + \mathbf{w}^i\|_2^2$ and $o_d^i:=\rho \|\mathbf{z}^{i}-\mathbf{z}^{i-1}\|_2^2$, accordingly. The iterations of ADMM can be terminated when both $o_p^i$ and $o_d^i$ are within $10^{-3}T\sqrt{N}$~\cite{Boyd10}. Figure~\ref{fig:dp1} demonstrates the cost convergence for \eqref{eq:EVS}, while Fig.~\ref{fig:dp2} the convergence of $o_p^i/T\sqrt{N}$ and $o_d^i/T\sqrt{N}$. As evidenced by Fig.~\ref{fig:dpncevs}, the global optimum is attained within 2,000 iterations.


\section{Conclusions}\label{sec:conclusions}
Given that optimal EV charging scales unfavorably with the fleet size and the number of control periods, decentralized charging protocols have been developed in this work.  A simple vehicle charging scheme has been devised based on the Frank-Wolfe iterations. This charging protocol exhibits provable $\mathcal{O}(\frac{1}{k})$ convergence, poses minimal computational requirements to EV controllers, enjoys privacy and security features, and attains a 100-times acceleration in terms of computational time over existing alternatives. To respect voltage and feeder transformer limits, network-constrained EV charging has been considered too. To achieve scalability, an ADMM-based solver has been built leveraging on an approximate grid model. The solver features closed-form updates and incorporates the scheduling protocol of vehicle charging. Numerical tests on real-world data verify the optimality and efficiency of the proposed decentralized schemes. Extensions to asynchronous ADMM and Frank-Wolfe updates, integrating the two optimization loops into a single update, and aiming for real-time EV scheduling schemes constitute current research directions.

\appendix\label{sec:appendix}
\begin{IEEEproof}[Proof of Prop.~\ref{pro:2}] The Lagrangian function of the convex problem in \eqref{eq:PQ0} reads
$L(\tilde{\mathbf{P}}_0, \tilde{\mathbf{Q}}_0,\nu)=\|\tilde{\mathbf{P}}_0-\breve{\mathbf{P}}_0^i\|_2^2 +\|\tilde{\mathbf{Q}}_0-\breve{\mathbf{Q}}_0^i\|_2^2+\nu\left[ (\mathbf{1}^{\top}\tilde{\mathbf{P}}_0)^2 + (\mathbf{1}^{\top}\tilde{\mathbf{Q}}_0)^2 - \bar{S}_f^2 \right]$. 
Because \eqref{eq:PQ0} satisfies Slater's condition (e.g., for $\tilde{\mathbf{P}}_0=\tilde{\mathbf{Q}}_0=\mathbf{0}$), strong duality holds~\cite{BoVa04}. If $(\tilde{\mathbf{P}}_0^*, \tilde{\mathbf{Q}}_0^*,\nu^*)$ are the optimal primal/dual variables, Lagrangian optimality yields:
\begin{subequations}\label{eq:kkt}
\begin{align}
\tilde{\mathbf{P}}_0^*-\breve{\mathbf{P}}_0^i+\nu^* \mathbf{1}\mathbf{1}^\top \tilde{\mathbf{P}}_0^*&=\mathbf{0}\label{eq:kkta}\\
\tilde{\mathbf{Q}}_0^*-\breve{\mathbf{Q}}_0^i+\nu^*\mathbf{1}\mathbf{1}^\top \tilde{\mathbf{Q}}_0&=\mathbf{0}.\label{eq:kktb}
\end{align}\
\end{subequations}
Premultiplying both sides of \eqref{eq:kkta}--\eqref{eq:kktb} by $\mathbf{1}^\top$ results in:
\begin{equation}\label{eq:rekkt}
\mathbf{1}^\top \tilde{\mathbf{P}}_0^*= \frac{\mathbf{1}^\top \breve{\mathbf{P}}_0^i}{1+3\nu^*}\quad \text{and} \quad \mathbf{1}^\top \tilde{\mathbf{Q}}_0^*= \frac{\mathbf{1}^\top \breve{\mathbf{Q}}_0^i}{1+3\nu^*}.
\end{equation}
Complementary slackness yields $\nu^*\big [ (\mathbf{1}^{\top}\tilde{\mathbf{P}}_0^*)^2 + (\mathbf{1}^{\top}\tilde{\mathbf{Q}}_0^*)^2 -  \bar{S}_f^2 \big]=0$, which from \eqref{eq:rekkt} and dual feasibility provides
\begin{equation}\label{eq:nu}
\nu^*=0~ \text{or}~\nu^*=\frac{1}{3}\left(\sqrt{(\mathbf{1}^\top \breve{\mathbf{P}}_0^i)^2+(\mathbf{1}^\top \breve{\mathbf{Q}}_0^i)^2}/\bar{S}_f - 1\right).
\end{equation}
The claim follows from primal feasibility, \eqref{eq:kkt}, and \eqref{eq:nu}.
\end{IEEEproof}

\bibliographystyle{IEEEtran}
\bibliography{myabrv,power}

\begin{thebibliography}{10}
\providecommand{\url}[1]{#1}
\csname url@samestyle\endcsname
\providecommand{\newblock}{\relax}
\providecommand{\bibinfo}[2]{#2}
\providecommand{\BIBentrySTDinterwordspacing}{\spaceskip=0pt\relax}
\providecommand{\BIBentryALTinterwordstretchfactor}{4}
\providecommand{\BIBentryALTinterwordspacing}{\spaceskip=\fontdimen2\font plus
\BIBentryALTinterwordstretchfactor\fontdimen3\font minus
  \fontdimen4\font\relax}
\providecommand{\BIBforeignlanguage}[2]{{%
\expandafter\ifx\csname l@#1\endcsname\relax
\typeout{** WARNING: IEEEtran.bst: No hyphenation pattern has been}%
\typeout{** loaded for the language `#1'. Using the pattern for}%
\typeout{** the default language instead.}%
\else
\language=\csname l@#1\endcsname
\fi
#2}}
\providecommand{\BIBdecl}{\relax}
\BIBdecl

\bibitem{RFK12}
P.~Richardson, D.~Flynn, and A.~Keane, ``Optimal charging of electric vehicles
  in low-voltage distribution systems,'' \emph{{IEEE} Trans. Power Syst.},
  vol.~27, no.~1, pp. 268--279, Feb. 2012.

\bibitem{YJCao12}
Y.~Cao, S.~Tang, C.~Li, P.~Zhang, Y.~Tan, Z.~Zhang, and J.~Li, ``An optimized
  {EV} charging model considering {TOU} price and {SOC} curve,'' \emph{{IEEE}
  Trans. Smart Grid}, vol.~3, no.~1, pp. 388--393, Mar. 2012.

\bibitem{SPR09}
S.~Shao, M.~Pipattanasomporn, and S.~Rahman, ``Challenges of {PHEV} penetration
  to the residential distribution network,'' in \emph{Proc. Power \& Energy
  Society General Meeting}, Calgary, Alberta Canada, Jul. 2009.

\bibitem{IP13}
Z.~Fan, ``A distributed demand response algorithm and its application to {PHEV}
  charging in smart grids,'' \emph{{IEEE} Trans. Smart Grid}, vol.~3, no.~3,
  pp. 1280--1290, Sep. 2012.

\bibitem{mah}
Z.~Ma, D.~Callaway, and I.~Hiskens, ``Decentralized charging control for large
  populations of plug-in electric vehicles,'' in \emph{Proc. Conf. on Decision
  and Control}, Atlanta, GA, Dec. 2010.

\bibitem{NGG}
N.~Gatsis and G.~B. Giannakis, ``Residential load control: {D}istributed
  scheduling and convergence with lost {AMI} messages,'' \emph{{IEEE} Trans.
  Smart Grid}, vol.~3, no.~2, pp. 770--786, Jun. 2012.

\bibitem{Rivera13}
J.~Rivera, P.~Wolfrum, S.~Hirche, C.~Goebel, and H.-A. Jacobsen, ``Alternating
  direction method of multipliers for decentralized electric vehicle charging
  control,'' in \emph{Proc. Conf. on Decision and Control}, Florence, Italy,
  Dec. 2013, pp. 6960--6965.

\bibitem{RQO14}
R.~Li, Q.~Wu, and S.~Oren, ``Distribution locational marginal pricing for
  optimal electric vehicle charging management,'' \emph{{IEEE} Trans. Power
  Syst.}, vol.~29, no.~1, pp. 203--211, Jan. 2014.

\bibitem{GTL11}
L.~Gan, U.~Topcu, and S.~Low, ``Optimal decentralized protocol for electric
  vehicle charging,'' \emph{{IEEE} Trans. Power Syst.}, vol.~28, no.~2, pp.
  940--951, May 2013.

\bibitem{karfopoulos2013multi}
E.~L. Karfopoulos and N.~D. Hatziargyriou, ``A multi-agent system for
  controlled charging of a large population of electric vehicles,''
  \emph{{IEEE} Trans. Power Syst.}, vol.~28, no.~2, pp. 1196--1204, May 2013.

\bibitem{NCEV}
J.~de~Hoog, T.~Alpcan, M.~Brazil, D.~Thomas, and I.~Mareels, ``Optimal charging
  of electric vehicles taking distribution network constraints into account,''
  \emph{{IEEE} Trans. Power Syst.}, vol.~30, no.~1, pp. 365--375, Jan. 2015.

\bibitem{MIPEV}
J.~Franco, M.~Rider, and R.~Romero, ``A mixed-integer linear programming model
  for the electric vehicle charging coordination problem in unbalanced
  electrical distribution systems,'' \emph{{IEEE} Trans. Smart Grid}, vol.~6,
  no.~5, pp. 2200--2210, Sep. 2015.

\bibitem{rev}
G.~Benetti, M.~Delfanti, T.~Facchinetti, D.~Falabretti, and M.~Merlo,
  ``Real-time modeling and control of electric vehicles charging processes,''
  \emph{IEEE Trans. Smart Grid}, vol.~6, no.~3, pp. 1375--1385, May 2015.

\bibitem{NCQ14}
N.~Chen, C.~W. Tan, and T.~Quek, ``Electric vehicle charging in smart grid:
  Optimality and valley-filling algorithms,'' \emph{{IEEE} J. Sel. Topics
  Signal Process.}, vol.~8, no.~6, pp. 1073--1083, Dec. 2014.

\bibitem{jaggi2013revisiting}
M.~Jaggi, ``Revisiting {Frank-Wolfe}: {P}rojection-free sparse convex
  optimization,'' in \emph{Proc. Intl. Conf. on Machine Learning}, Atlanta, GA,
  Jun. 2013, pp. 427--435.

\bibitem{BoVa04}
S.~Boyd and L.~Vandenberghe, \emph{Convex Optimization}.\hskip 1em plus 0.5em
  minus 0.4em\relax New York, NY: Cambridge University Press, 2004.

\bibitem{donald1999art}
D.~E. Knuth, \emph{The Art of Computer Programming: {Sorting and
  Searching}}.\hskip 1em plus 0.5em minus 0.4em\relax Boston, MA: Pearson
  Education, 1998, vol.~3.

\bibitem{Be15}
D.~P. Bertsekas, \emph{Convex Optimization Algorithms}.\hskip 1em plus 0.5em
  minus 0.4em\relax Belmont, MA: Athena Scientific, 2015.

\bibitem{LAGan}
\BIBentryALTinterwordspacing
L.~Gan and S.~H. Low. (2014) Convex relaxations and linear approximations for
  optimal power flow in multiphase radial networks. [Online]. Available:
  \url{http://arxiv.org/pdf/1406.3054v1.pdf}
\BIBentrySTDinterwordspacing

\bibitem{HaozhuD}
E.~Dall'Anese, H.~Zhu, and G.~B. Giannakis, ``Distributed optimal power flow
  for smart microgrids,'' \emph{{IEEE} Trans. Smart Grid}, vol.~4, no.~3, pp.
  1464--1475, Sep. 2013.

\bibitem{VZ}
V.~Kekatos, L.~Zhang, G.~B. Giannakis, and R.~Baldick, ``Voltage regulation
  algorithms for multiphase power distribution grids,'' \emph{{IEEE} Trans.
  Power Syst.}, vol.~PP, no.~99, pp. 1--11, 2015.

\bibitem{Boyd10}
S.~Boyd, N.~Parikh, E.~Chu, B.~Peleato, and J.~Eckstein, ``Distributed
  optimization and statistical learning via the alternating direction method of
  multipliers,'' \emph{{F}ound. {T}rends {M}ach {L}earning}, vol.~3, pp.
  1--122, 2010.

\bibitem{KeGi12}
V.~Kekatos and G.~B. Giannakis, ``Distributed robust power system state
  estimation,'' \emph{{IEEE} Trans. Power Syst.}, vol.~28, no.~2, pp.
  1617--1626, May 2013.

\bibitem{QiuyuPeng}
Q.~Peng and S.~Low, ``Distributed algorithm for optimal power flow on a radial
  network,'' in \emph{Proc. Conf. on Decision and Control}, Venice, Italy, Dec.
  2014, pp. 167--172.

\bibitem{Federalhighway}
\BIBentryALTinterwordspacing
Federal highway administration. US Department of Transportation. [Online].
  Available: \url{http://nhts.ornl.gov/2009/ pub/stt.pdf}
\BIBentrySTDinterwordspacing

\bibitem{EVload}
K.~Qian, C.~Zhou, M.~Allan, , and Y.~Yuan, ``Modeling of load demand due to
  {EV} battery charging in distribution systems,'' \emph{{IEEE} Trans. Power
  Syst.}, vol.~26, no.~2, pp. 802--810, May 2011.

\bibitem{TestFeeder}
\BIBentryALTinterwordspacing
Distribution test feeders. IEEE Power \& Energy Society. [Online]. Available:
  \url{http://nhts.ornl.gov/2009/ pub/stt.pdf}
\BIBentrySTDinterwordspacing

\end{thebibliography}
\end{document}